\documentclass{article}

\usepackage{amsfonts}
\usepackage{amssymb}
\usepackage{latexsym}
\usepackage{ifthen}

\newtheorem{thm}{Theorem}
\newtheorem{lem}{Lemma}

\newtheorem{cor}{Corollary}
\newcommand{\ra}{\rightarrow}
\newcommand{\tx}{\textrm}

\begin{document}

\title{An elementary purely algebraic approach to generalized functions and operational calculus}
\author{Vakhtang Lomadze}
\date{~~ Andrea Razmadze Mathematical Institute, Mathematics Department of I. Javakhishvili Tbilisi State University,  Georgia}
\maketitle

{\bf Abstract.}\
The space  of Schwartz distributions of finite order is represented as a factor space
of the space  of, what we call, Mikusinski functions. The point of Mikusinski functions is that they admit a  multiplication by
convergent Laurent series. It is shown that this multiplication provides a natural
simple basis for Heaviside's  operational calculus.


\bigskip


\section{Introduction}

  In his  seminal works in electromagnetic theory, O. Heaviside  developed formal rules for dealing with the differentiation operator, much of which he arrived at intuitively.    A mathematical basis for his operational calculus was done by mathematicians later with the aid of the Laplace transform.
 A different explanation of Heaviside's calculus, based upon
 the convolution ring of continuous functions, was presented by Mikusinski \cite{m}.
From this ring, which is commutative and without zero divizors, Mikusinski goes on to define the field of fractions. These fractions he  calls operators.

 In this note, we modify Mikusinski's approach by starting with the module structure  of continuous functions over  the ring of convergent formal series (say, in $t$). This module  is  {\em torsion free}.  Imitating Mikusinski, we  apply the fraction space construction and  obtain this way a linear space over the field of convergent Laurent series. Elements of this space, which we call Mikusinski functions, constitute  a very small part of
the
field of Mikusinski operators. Viewing the space of these functions as a module over the polynomial ring
(in $s=t^{-1}$)
and taking the factor module by "derivatives" of constant functions, we get a module  that is canonically isomorphic
to the module ${\cal D}^\prime_{fin}$, the module of Schwartz distributions of finite order.

The goal of this note is to demonstrate  that with such a representation of ${\cal D}^\prime_{fin}$ we have a natural, simple and rigorous theory for Heaviside's operational calculus.

Throughout, $I$ is a fixed interval of real axis containing $0$ and  $t$  an indeterminate. The symbol $\bf{1}$ will stand for the unit function on $I$.
Let $C(I)$ be the space of all complex-valued continuous  functions defined on the interval $I$ and let $\mathbb{C}\{t\}$ be
the ring of complex coefficient  convergent formal
series in $t$.
(We remind that a formal series $\sum_{i\geq 0}b_it^i$ is said to be convergent if $\sum_{i\geq 0}|b_i|\varepsilon^i < +\infty$
 for some positive real $\varepsilon$.) Let $\mathbb{C}(\{t\})$ denote the fraction field of $\mathbb{C}\{t\}$, the field of convergent Laurent series.
It is worth mentioning that $\mathbb{C}(\{t\})$ contains as a subring the polynomial ring $\mathbb{C}[s]$, where $s=t^{-1}$.

For every continuous  function $u\in C(I)$, let $J(u)$ denote the function defined by
$$
J(u)(x)=\int_0^xu(\alpha)d\alpha,\ \ \ x\in I.
$$

\section{Mikusinski functions}

If
 $g=\sum_{i\geq 0}b_it^i$ is a convergent formal series and $u$ is a continuous function on $I$, then the series
$$
\sum_{n\geq 0} b_nJ^n (u)
$$
 converges uniformly on every compact subinterval of $I$ that contains $0$. We therefore can  define
 the product $gu$ by setting
$$
gu=\sum_{n\geq 0} b_nJ^n (u).
$$

This multiplication   makes $C(I)$ a {\em module} over $\mathbb{C}\{t\}$.
 (A module is a "linear space" over a ring!) It is easily seen that this module is torsion free. Indeed, suppose that $gu=0$, where $g\neq 0$. Write $g=t^ng_0$, where $g_0$ has nonzero
free coefficient. By the theorem on units (see \cite{r}), $g_0$ is invertible. Multiplying $gu=0$ by $g_0^{-1}$, we get
that $t^nu=0$. Because $J$ is an injective operator, it follows that $u=0$.

{\bf Definition}. We  define the Mikusinski space $M(I)$
 as the fraction space of the module $C(I)$.

  By definition,  $M(I)$ is a linear space over $\mathbb{C}(\{t\})$, and every Mikusinski  function is represented by a ratio $ u/g$ with $u\in C(I)$ and $g\in \mathbb{C}\{t\},\ g\neq 0$.
Two such ratios $u_1/g_1$ and $u_2/g_2$ represent the same Mikusinski function if
$$
g_2u_1=g_1u_2.
$$

 Because $C(I)$ is torsion free, the canonical map $u \mapsto u/1$ is injective.  Therefore we can make the identification
$$
u=u/1.
$$
With this identification,
 a Mikusinski  function $w$  can be written in the form
$$w=s^nu,$$
where $u\in C(I)$.
 Indeed, if $w=u/g$, then, writting $g=t^ng_0$ with $g_0$ as above, we have
$$\frac{u}{g}=s^n\frac{g_0^{-1}u}{1}=s^n(g_0^{-1}u).$$
We therefore have:
$$
C(I)\subseteq sC(I)\subseteq s^2C(I)\subseteq s^3C(I)\subseteq \cdots  \ \ \tx{and}\ \ M(I)=\cup_{n\geq 0} s^nC(I).
$$

{\em Remark}. If one wants, one can define the Mikusinski space $M(I)$ as the inductive limit
of the sequence
$$
C(I)\stackrel{J}{\ra} C(I) \stackrel{J}{\ra}C(I) \stackrel{J}{\ra} C(I)\stackrel{J}{\ra}  \cdots .
$$

Extend the operator $J$ to Mikusinski functions by setting
$$
J(w)=tw.
$$
This way, the integration operator becomes bijective.
The following theorem says that the  space $M(I)$ is  the smallest  extension of $C(I)$
   in which every  function is an  integral.
\begin{thm} Every continuous function is an $n$-fold integral of a Mikusinski function.
\end{thm}
{\em Proof}. This is obvious. Indeed, if $u\in C(I)$, then
$$
u=t^n(s^nu)=J^n(s^nu).
$$
$\quad\Box$

\section{Schwartz distributions (of finite order)}

Since multiplication by $t$ is the integration operator, a natural idea is  to interpret multiplication by $s=t^{-1}$ as the differentiation operator. Postulating that constant functions must have zero derivatives,
we are led to the following definition.

{\bf Definition}. Define the generalized function space $G(I)$ to be
$$
G(I)\ =\ M(I)/N(I),
$$
where $N(I)=\{f{\bf 1}\ |\ f\in s\mathbb{C}[s]\}$.

Since $N(I)$ is a $\mathbb{C}[s]$-submodule of $M(I)$, $G(I)$ has a structure of
a module over $\mathbb{C}[s]$.
 In particular, we can multiply generalized functions by $s$.

{\bf Definition}. Define the differential operator $D: G(I)\ra G(I)$ by setting
$$
D(\xi)=s\xi,\ \ \ \xi\in G(I).
$$

We need the following lemma.
\begin{lem} \
$
C(I) \cap N(I) = \{0\}.
$
\end{lem}
{\em Proof}. Assume there is $u\in C(I)$ that is not zero and belongs to $N(I)$. We
then have
$$
u = (a_0s^n + \cdots + a_{n-1}s){\bf 1}
$$
with $n\geq 1$, $a_i\in \mathbb{C}$ and $a_0\neq 0$. Multiplying  both of the sides by $t^{n+1}$, we get
$$
t^{n+1}u = (a_0t+ \cdots + a_{n-1}t^n){\bf 1}.
$$
It follows  that
$$
t(t^nu - (a_1+ \cdots + a_{n-1}t^{n-1}){\bf 1}) = a_0{\bf 1}.
$$
On the left, we have a continuous function having value 0 at 0. Hence, $a_0=0$, which is a contradiction.
$\quad\Box$

It follows from the previous lemma that the canonical map $C(I)\ra G(I)$ defined by
$$
u \mapsto u\ (\tx{mod}\ N(I))
$$
is injective. This permits us to  identify $C(I)$ with its image under this map. In other words, for every
$u\in C(I)$, we can make the identification
$$
u= u\ (\tx{mod}\ N(I)).
$$

Notice that the derivative $D(u)$ of a classical  continuously differentiable function $u$ coincides with the ordinary derivative $u^\prime$.
 Indeed,  by the Newton-Leibniz formula
$u=Ju^\prime + u(0){\bf 1}$,
and consequently
$$
su\ (\tx{mod}\ N(I)) =(u^\prime +su(0){\bf 1})\ (\tx{mod}\ N(I)) =u^\prime (\tx{mod}\ N(I))=u^\prime.
$$

The following theorem says that the space $G(I)$ is  the smallest extension of $C(I)$ in which every function has a derivative.

\begin{thm}   Any generalized function is an $n$-fold derivative of a continuous function.
\end{thm}
{\em Proof}. This is obvious. Indeed, for every continuous function $u$ and every nonnegative integer $n$, we have
$$
s^nu\ (\tx{mod}\ N(I))=D^n (u\ (\tx{mod}\ N(I)))= D^nu.
$$ $\quad\Box$

As an immediate consequence, we have the following corollary.
\begin{cor}
 $G(I)$ is canonically isomorphic to ${\cal D}^\prime_{fin}(I)$, the space of Schwartz distributions of finite order.
\end{cor}

This corollary permits us to interpret generalized functions as Schwartz distributions of finite order.

{\em Remark}. We remind that, basically, distributions have been introduced by Schwartz  in order to be able to differentiate all continuous functions
(see page 72 in Schwartz \cite{s}).

 We close the section with  the following important theorem.

\begin{thm}  Let $f\in \mathbb{C}[s]$ be a polynomial. The equation
$$
f(D)\xi=\omega
$$
has a solution $\xi\in G(I)$ for every $\omega\in G(I)$.
\end{thm}
{\em Proof}. This is obvious. Indeed,   $M(I)$ is a divisible  $\mathbb{C}[s]$-module, and consequently $G(I)$ also is  divisible.
To be more precise, if $\omega$ is represented by a Mikusinski function $w$, then $$\xi=\frac{1}{f}w \ \tx{mod} N(I)$$ is a solution.
$\quad\Box$

\section{Linear constant coefficient differential equations}

We already know how to find a particular solution of a nonhomogeneous linear constant coefficient differential equation.
So, we only need to consider the homogeneous case.

Given $g\in \mathbb{C}\{t\}$, define  an entire analytic function $E(g)$ by setting $$E(g)=g{\bf 1}.$$

\begin{thm} Let $f\in \mathbb{C}[s]$ be a polynomial of degree $d$. The solutions of the linear differential equation
$$
f(D)\xi=0,\ \ \ \ \xi\in G(I)
$$
 are given by the formula
$$
\xi=E(\frac{sr}{f}),
$$
where $r$ runs over the polynomials in $\mathbb{C}[s]$ of degree $\leq d-1$.
\end{thm}
{\em Proof}. It is  easily seen
 that all these functions are solutions. Indeed,
$$
f(D)\xi=f\frac{sr}{f}{\bf 1}=sr {\bf 1}=0.
$$
To show the converse, assume that $\xi$ is  a solution of our equation and assume that a Mikusinski function $s^nx$ represents it.  Then
$$
fs^nx\in N(I).
$$
This means that there exist complex numbers  $a_0, a_1, \ldots , a_k$ such that
$$
fs^nx=s(a_0+a_1s+ \cdots + a_ks^k){\bf 1}.
$$
Multiplying both sides by $t^{d+n}$, from this we get
$$
(ft^d)x=(a_0t^{d+n-1}+\dots +a_{d+n-1}){\bf 1} + (a_{d+n}s+\cdots  ){\bf 1}.
$$
Since the left side is an ordinary  continuous function, we can see that $a_k=0$ for all $k\geq d+n$.
We therefore have
$$
s^nx=\frac{s(a_0+a_1s+ \cdots + a_{d+n-1}s^{d+n-1})}{f}{\bf 1}.
$$
By the Euclidean division theorem, there exists a polynomial $r$ of degree $\leq d-1$ such that
$$
\frac{a_0+a_1s+ \cdots + a_{d+n-1}s^{d+n-1}}{f} \equiv \frac{r}{f} \ \tx{mod}\ \mathbb{C}[s].
$$
So that
$$
\xi=E(\frac{sr}{f}).
$$
The proof is complete. $\quad\Box$

\end{document}